\documentclass[sn-mathphys,Numbered]{sn-jnl}
\usepackage{graphicx}%
\usepackage{multirow}%
\usepackage{amsmath,amssymb,amsfonts}%
\usepackage{amsthm}%
\usepackage{mathrsfs}%
\usepackage[title]{appendix}%
\usepackage{xcolor}%
\usepackage{textcomp}%
\usepackage{manyfoot}%
\usepackage{booktabs}%
\usepackage{algorithm}%
\usepackage{algorithmicx}%
\usepackage{algpseudocode}%
\usepackage{listings}%
\usepackage{natbib}
\usepackage{bm}
\usepackage{caption}

\newtheorem{theorem}{Theorem}

\theoremstyle{thmstyletwo}%
\newtheorem{remark}{Remark}%

\theoremstyle{thmstylethree}%
\begin{document}
	
	\title[Article Title]{Uniform boundary observability for the spectral collocation of the linear elasticity system}

	\author*[1]{\fnm{Somia} \sur{Boumimez}}\email{soumia.boumimez@gmail.com}
	\author[1]{\fnm{Carlos} \sur{Castro}}\email{carlos.castro@upm.es}
	\affil*[1]{M2ASAI Universidad Polit\'ecnica de Madrid, Departamento de Matem\'atica e Inform\'atica, ETSI Caminos, Canales y Puertos, 28040 Madrid, Spain}
	
	\abstract{	A well-known boundary observability inequality for the elasticity system establishes that the energy of the system can be estimated from the solution on a sufficiently large part of the boundary for sufficiently large time. 
		This inequality is relevant in different contexts as the exact boundary controllability, the boundary stabilization or some inverse source problems. Here we show that a corresponding  boundary observability inequality for the spectral collocation approximation of the linear elasticity system in a $d$-dimensional cube also holds, uniformly with respect to the discretization parameter. This property is essential  to prove that natural numerical approaches of the previous problems based on replacing the elasticity system by the collocation discretization will give successful approximations of the continuous counterparts. As an application we obtain the boundary controllability of the discrete system resulting when approximating the elasticity system with this numerical method, uniformly with respect to the discretization parameter. We also give numerical evidences of the convergence of these discrete controls to a boundary control of the limit $ 2d-$elasticity system in a square domain.}
	\keywords{Observability, Controllability, Elasticity systems, Spectral collocation methods.}
	\maketitle
	
	\section{Introduction}
		\hspace*{0.5cm}We consider the free vibrations of a $d$-dimensional homogeneous and isotropic elastic body occupying a bounded domain $ \Omega =(-1,1)^d \subset \mathbb{R}^{d},(d\ge2) $ with boundary $ \partial \Omega,$
	\begin{align}\label{adjcont}
		\begin{cases}
			\left(\bm{\phi}_{tt}- \Delta^*\bm{\phi}\right) (t,{\bf x})={\bf 0}& (t,{\bf x}) \in (0,T)\times\Omega, \\
			\bm {\phi}(t,{\bf x})={\bf 0},& (t,{\bf x}) \in  (0,T)\times  \partial\Omega\\
			( \bm {\phi}(0,{\bf x}),\bm {\phi}_{t}(0,{\bf x}))=(\bm {\phi^0}({\bf x}),\bm{\phi^1}({\bf x})),\hspace*{0.5cm}&{\bf x} \in  \Omega,
		\end{cases}
	\end{align}
	where  $ \bm{\phi}(t,{\bf x})$ is the displacement of the material point $ {\bf x} $ at time $ t $, $ \Delta^*=\mu\Delta+(\lambda+\mu)\nabla \mbox{div}$,  $ \lambda, \mu > 0  $ are Lam$\acute{e} $  parameters and $ (\bm {\phi^0,\phi^1})\in (H_0^1(\Omega))^d\times (L^2(\Omega))^d $ are the initial displacement and velocity in the usual energy space. It is  well known that problem \eqref{adjcont} is observable from a sufficiently large part of the boundary $\Gamma\subset \partial \Omega$ and for sufficiently large time $T$  (see \cite{cit6} and \cite{cit44b}). More precisely,  there exists a constant $ C>0 $ such that 
	\begin{equation}		\label{obs_ineq}
		C\left|(\bm{\phi}^{\bf 0},\bm{\phi}^{\bf 1})\right|^{2}_{(H_0^1)^d\times (L^2)^d}\le\int_{0}^{T}\int_{\Gamma}\left|\mu \dfrac{\partial\bm{\phi}}{\partial{\bm \nu}}+(\lambda+\mu){\bm \nu} \mbox{div} \bm{\phi}\right|^{2} d\gamma dt,
	\end{equation}
	for any initial data  $ (\bm {\phi^0,\phi^1})$, where $ {\bm \nu} $ is the unit normal vector to $ \Gamma $, directed towards the exterior of $ \Omega $.	
	Although this is true in more general domains and situations here we focus on the particular case  $T>\dfrac{4\sqrt{d}}{\sqrt{\mu}}$ and $\Gamma$  the part of the boundary constituted by $d$ adjacent non-opposite faces:  
	$
	\Gamma= \left\{ (x_1,..,x_d),\quad s.t. \quad x_j=1, \; \mbox{ for at least one } j=1,...,d \right\}.
	$ 
	The fact that we only consider rectangular domains is due to the particular spectral method that we analyze below. On the other hand, the choice of $\Gamma$ considered here is somehow optimal in the sense that  the observability inequality fails for any smaller open subset of $\Gamma$. Inequality \eqref{obs_ineq} is essential to establish boundary controllability properties, stabilization or inverse source problems for the elasticity system from boundary observations (see \cite{cit6}, \cite{cit44b} and \cite{T}, for inverse source problems in simpler models).  
	
	Regarding the numerical approximation of such problems the natural approach consists in replacing the continuous problem with a suitable discrete approximation. In particular, this requires a finite dimensional discrete approximation of the elasticity system \eqref{adjcont} and a corresponding observability inequality for the discrete system. Moreover, this inequality must be uniform with respect to the discretization parameter in order to establish any convergent result. It turns out that such discrete observability inequalities hold but are not uniform for the usual discretization coming from finite elements and finite differences. This has been extensively analyzed for the simpler scalar wave equation (see the review papers in \cite{cit7}, \cite{EZ} and \cite{cit11} together with the references therein). 
	
	Here we propose a spectral collocation method based on a polynomial approximation. For the scalar wave equation in one and two dimensions this method provides a uniform observability inequality, as long as we introduce an extra term (see \cite{Som}). We also refer to \cite{urqita} where the authors show that without this extra term, the inequality cannot be uniform. We show that a similar uniform observability inequality holds for the discrete elasticity system in any dimension. As far as we know, this is one of the first results in this context valid for a convergent discretization of the elasticity system. For a 
	similar result using a less classical approach based on suitable stabilized space-time finite elements we refer to the recent work \cite{BM}.

	As an application we show how this result provides a method to construct numerical approximations of the boundary control for the elasticity system in dimension $d=2$. This problem was also addressed in \cite{BM} using the space-time finite element approach mentioned above. It is also worth mentioning reference \cite{RF} where the authors find numerical approximations of these boundary controls, in some particular cases, using the energy disipation of the solutions of the elasticity system in a extended domain. One of the main advantages of our approach is the fact that the method recovers the high accuracy inherent to the spectral approximation.  
	
	The rest of this note is organized as follows. In Section 2, we introduce the numerical method and prove the associated uniform observability inequality. In Section 3 we show how to apply this result to obtain convergent approximations of the boundary control for the elasticity system. 
	
\section{Spectral collocation method}
\hspace*{0.5cm} For the background and details on this method to approximate the solutions of the elasticity system we refer the reader to \cite{ex5}.
	Let $ {\bf N}= (N,...,N) \in \mathbb{ N}^d $  be the number of collocation points in each variable $ x_j  $ that we assume to be the same. We can consider more general situations with different number of points in each dimension but this is not relevant in our analysis and would make the notation more involved. We also consider $ C=\{{\bf P_ i}=(x^{k_1}_{1},...,x^{k_d}_{d}),(0,...,0)\le {\bf i}=(k_1,...,k_d)\le (N,...,N)\} $ the Legendre-Gauss-Lobatto (LGL) nodes in $ \Omega $ that are the roots of the polynomial
	$\prod_{j=1}^{d}(1-x_{j}^{2})\partial_{x_{j}}L_{N}(x_{j}),$ 
	where $ L_{k} $ is the k-th Legendre polynomial in $ (-1, 1)  $ (e.g. \cite{cit9}). We divide $C$  into interior and boundary nodes, i.e. $ C = C^{\Omega}\cup C^{\partial\Omega}$ where
	$C^{\Omega}=C\cap \Omega =\{ {\bf P_i, i\in I}_{\Omega}\} $ and
	$ C^{\partial\Omega}=C\cap \partial\Omega=\{ {\bf P_i, i\in I}_{\partial\Omega}\}$, and ${\bf I}_{\Omega},\; {\bf I}_{\partial\Omega}$ are the sets of indexes corresponding to the interior and boundary collocation nodes respectively. We denote $ {\bf I}={\bf I}_{\Omega}\cup {\bf I}_{\partial \Omega} $.
	
	Let $ \mathbb{P}_{\bf N}(\Omega)  $ be the space of  polynomials of degree
	at most $ N $ in the $ x_{j} $-variable, $ j=1,...,d $ and let $ \mathbb{P}_{\bf  N}^{Di}(\Omega)  $ be the subspace of $ \mathbb{P}_{\bf N}(\Omega)  $ of those vanishing on the boundary $ \partial \Omega $. Consider the following collocation approximation of system \eqref{adjcont}: Find $\bm{\phi}^{\bf N} \in (C^\infty([0,T];\mathbb{P}_{\bf N}(\Omega)))^d$ such that
	\begin{align}\label{adj-dis}
		\begin{cases}
			(\bm{\phi}^{\bf N}_{tt}-\Delta^* \bm{\phi}^{\bf N}) (t,{\bf P_{\bf i}})={\bf 0}, & (t,{\bf P_{\bf i}}), \in (0,T)\times C^{\Omega},\\
			\bm{\phi}^{\bf N}(t,{\bf P_{\bf i}})={\bf 0}, &  (t,{\bf P_{\bf i}}),\in (0,T)\times C^{\partial\Omega},\\
			\bm{\phi}^{\bf N}(0,\cdot)=\bm{\phi}^{\bf 0, N},\quad \bm{\phi}_{t}^{\bf N}(0,\cdot)=\bm{\phi}^{\bf 1,N},&
		\end{cases}
	\end{align}
	where $ (\bm{\phi}^{\bf 0,N},\bm{\phi}^{\bf 1,N})\in (\mathbb{P}^{Di}_{\bf N}(\Omega))^d\times (\mathbb{P}^{Di}_{\bf N}(\Omega))^d$. Note that \eqref{adj-dis} is a second order system of ODE with $d(N+1)^d$ equations and unknowns, namely the coefficients of the polynomial $\bm{\phi}^{\bf N}$.  
	The main result in this paper is the following uniform discrete version of \eqref{obs_ineq}.	
	
\begin{theorem}\label{obser} System \eqref{adj-dis} is uniformly observable from the boundary $\Gamma$ in time $ T >\dfrac{4\sqrt{d}(2+N^{-1})^d}{\sqrt{\mu}}  $. More precisely, there exist two constants $  C_1,C_2 > 0 $, independent of $ {\bf N} $, such that 
		\begin{align}\label{obser-inq}
			\begin{split}
				C_1\left|(\bm{\phi}^{\bf 0,N},\bm{\phi}^{\bf 1, N})\right|^{2}_{(H_0^1)^d\times (L^2)^d}&\le\int_{0}^{T}\int_{\Gamma}\left|\mu \dfrac{\partial\bm{\phi}^{\bf N}}{\partial{\bm \nu}}+(\lambda+\mu){\bm \nu} \mbox{div} \bm{\phi}^{\bf N}\right|^{2} d\gamma dt\\
				&+C_{2}\int_{0}^{T}\int_{\partial \Omega}\left|\mu \dfrac{\partial^{2}\bm{\phi}^{\bf N}}{\partial{\bm \nu}^{2}}+(\lambda+\mu){\bm \nu}\dfrac{\partial\mbox{div}{\bm\phi}^{\bf N}}{\partial{\bm \nu}}\right|^{2} d\gamma dt,
			\end{split}
		\end{align}
		for  all initial data  $ (\bm{\phi}^{\bf 0,N},\bm{\phi}^{\bf 1,N})\in (\mathbb{P}^{Di}_{\bf N}(\Omega))^d\times (\mathbb{P}^{Di}_{\bf N}(\Omega))^d$.
	\end{theorem}
	\begin{remark}\label{remark2}
		Note that the observability inequality \eqref{obser-inq} is essentially the same as in \eqref{obs_ineq} with an extra term on the right hand side. It is easy to see that, at least formally, this extra term vanishes for smooth solutions of system \eqref{adjcont}. We can simply write the first equation in the elasticity system \eqref{adjcont} in local coordinates at the boundary and use the boundary conditions. This is not the case for the solutions of \eqref{adj-dis} where the elasticity system is prescribed only at the interior collocation points. This extra term must be taken into account when using this estimate in the applications. 
	  \end{remark}
  
	Let us introduce some notation before giving the proof.	We define the following
	discrete inner product that approximates the $ L^{2}(\Omega) $ one:
	\begin{align}
		(w,z)_{\bf N}=\sum_{{\bf i \in I}}(wz)({\bf P_ i})\; \omega_{\bf i}, \quad w,z\in \mathbb{P}_{\bf N}(\Omega).
	\end{align}
	Here  $ \omega_{\bf i}=\prod_{j=1}^{d}\omega_{k_j}$ where $ \omega_{k_j} $ is the discrete weight associated with the 1-d Legendre-Gauss-Lobato (LGL) quadrature formula  (e.g.\cite{ex5}, Chap. 2). Owing to the exactness of this quadrature,
	\begin{align}\label{intn}
		(w,z)_{\bf N}=\int_{\Omega} wz\ d{\bf x} \ \mbox{for\ all}\ w,z\ \mbox{such\ that}\ wz\in\mathbb{P}_{\bf 2 N-1}(\Omega), \quad {\bf x}=(x_1,...,x_d).
	\end{align}
	Moreover, the discrete norm $ \left\|\cdot\right\| _{\bf N} =\sqrt{(z,z)_{\bf N}} $ is uniformly equivalent to the $ \left|\cdot\right|_{L^{2}} -$norm  in $ \mathbb{P}_{\bf N}(\Omega) $ 
	(\cite{ex5}, Chapter 9). In fact, for the constants $ C_{1}=1 $, and $ C _{2}=(2+N^{-1})^d$, 
	\begin{align}\label{equiv-n}
		C_{1}\left|p\right|^{2}_{L^{2}}\le\left\|p\right\|^{2}_{\bf N}\le C _{2}\left|p\right|^{2}_{L^{2}}, \ \forall p\in \mathbb{P}_{\bf N}(\Omega) .	
	\end{align}
	We denote by $ \Psi_{k_j}(x_j),j=1,...,d,\; k_j=0,...,N $ the Lagrange polynomial which is $ 1 $ at $x_{j}^{k_j} $ and $ 0 $  at all the other collocation points. Observe that $\{\Psi_{\bf i}({\bf x})=\prod_{j=1}^{d}\Psi_{k_j}(x_j),\; {\bf i}=(k_1,...,k_d)\in {\bf I}\}$ constitutes a basis in $ \mathbb{P}_{\bf N}^{Di}(\Omega)  $.
	
	\begin{proof}[\textbf{Proof of Theorem \ref{obser}.}] The main idea is to observe that the solutions of  \eqref{adj-dis} solve the following equivalent continuous system, 
		\begin{align}\label{adj-equiv}
			\begin{cases}
				(\bm{\phi}^{\bf N}_{tt}-\Delta^*\bm{\phi}^{\bf N})(t,{\bf x}) =-\sum_{{\bf i\in  I}_{\partial\Omega}}\Delta^*{\bm\phi}^{\bf N}(t,{\bf P_i}) \Psi_{\bf i}({\bf x}), & (t,{\bf x})\in (0,T)\times\Omega, \\
				\bm{\phi}^{\bf N}(t,{\bf x})={\bf 0},& (t,{\bf x}) \in (0,T)\times\partial\Omega,\\
				\bm{\phi}^{\bf N}(0,{\bf x})=\bm{\phi}^{\bf0, N},\bm{\phi}_{t}^{\bf N}(0,{\bf x})=\bm{\phi}^{\bf 1, N}, & {\bf x} \in \Omega.
			\end{cases}
		\end{align}
		In fact, this is easily seen by writing $ \bm{\phi}^{\bf N} $ and $\bm{\phi}^{\bf N}_{tt}-\Delta^*\bm{\phi}^{\bf N}$ in the Lagrangian basis, since they are polynomial of degree $ N $ in each variable, and using system \eqref{adj-dis}.
		Note that system \eqref{adj-equiv} is a perturbation of the original system  \eqref{adjcont} so that we can adapt the continuous multipliers proof in Lions \cite{cit44b}, estimating the extra nonhomogeneous right hand side in \eqref{adj-equiv}. Given ${\bf x^0}=(-1,...,-1) \in \mathbb{R}^{d}$ and ${\bf m}({\bf x})\in \mathbb{R}^d$ with components $ m_j({\bf x})=x_{j}- x^0_{j},\;j=1,...,d $ it is clear that for all $ {\bf x} \in \partial\Omega\backslash\Gamma: $ ${\bf m} \bm{\cdot}{\bm \nu}={\bf 0} $. 	Let us set 	$
		X= \int_{\Omega}\bm{\phi}^{\bf N}_{t}\bm{\cdot} m_{j} \dfrac{\partial\bm{\phi}^{\bf N}}{\partial x_{j}} d{\bf x}\bigg| ^{T}_{0}$, $
		Y=\int_{\Omega} \bm{\phi}^{\bf N}_{t}\bm{\cdot}\bm{\phi}^{\bf N}d{\bf x}\bigg| ^{T}_{0}.$
		Multiplying scalarly the first vector equation in \eqref{adj-equiv} by the vector $ m_j \dfrac{\partial\bm{\phi}^{\bf N}}{\partial x_{j}} $ (here the repeated index $j$ stands for the sum in $j=1,...,d$) and integrating by parts, the left hand side can be simplified as follows, 
		\begin{align*}
			& X+\dfrac{d-1}{2}Y-\int_{0}^{T}\int_{\partial\Omega} \dfrac{\bm {m\cdot\nu}}{2}\bigg( \mu\left| \dfrac{\partial \bm{\phi}^{\bf N}}{\partial {\bm\nu}}\right| ^{2}+(\lambda+\mu)\left| \mbox{div} \bm{\phi}^{\bf N} \right| ^{2}\bigg)d\gamma dt\\
			&+\dfrac{1}{2}\int_{0}^{T} \int_{\Omega}(\left| \bm{\phi}_{t}^{\bf N}\right| ^{2}+\mu\left| \nabla\bm{\phi}^{\bf N}\right| ^{2}+(\lambda+\mu)\left| \mbox{div}\bm{\phi}^{\bf N}\right| ^{2}) d{\bf x} dt.
		\end{align*}
		Using  the fact that $ \bm {m \cdot \nu}=0 $ on $\ \partial\Omega\backslash\Gamma $ and $sup_{{\bf x}\in \Gamma}\ {m_j }{ \nu_j}=2\sqrt{d}$, we obtain
		\begin{align}\label{lions1}
			\begin{split}
				&\dfrac{1}{2}\int_{0}^{T} \int_{\Omega}\left(\left| \bm{\phi}_{t}^{\bf N}\right| ^{2}+\mu\left| \nabla\bm{\phi}^{\bf N}\right| ^{2}+(\lambda+\mu)\left| \mbox{div}\bm{\phi}^{\bf N}\right| ^{2}\right) d{\bf x} dt\le\left| X+\dfrac{d-1}{2}Y\right|\\
				&+\dfrac{2\sqrt{d}}{2}\int_{0}^{T}\int_{\Gamma}\bigg( \mu\left| \dfrac{\partial \bm{\phi}^{\bf N}}{\partial \bm{\nu}}\right| ^{2}+(\lambda+\mu)\left| \mbox{div} \bm{\phi}^{\bf N} \right| ^{2}\bigg)d\gamma dt\\
				&+\left|\int_{0}^{T}\int_{\Omega}\left( \sum_{{\bf i\in  I}_{\partial\Omega}}\bigg(\mu\dfrac{\partial^{2} \bm{\phi}^{\bf N} }{\partial \bm{\nu}^{2}}+(\lambda+\mu)\bm{\nu}\dfrac{\partial \mbox{div}\bm{\phi}^{\bf N} }{\partial \bm{\nu}}\bigg)(t,{\bf P_{\bf i}})\Psi_{\bf i}(\bf x)\right) \bm{\cdot}\left(m_{j}\dfrac{\partial\bm{\phi}^{\bf N}}{\partial x_{j}}\right)d{\bf x} dt\right|.
			\end{split}
		\end{align}
		We now estimate each one of the terms in this expression. We start with the left hand side in \eqref{lions1}. define the discrete energy 
		\begin{align}\label{disenergy}
			E^{\bf N}(t)=\dfrac{1}{2}\bigg( \left\| \bm{\phi}^{\bf N}_{t}\right\|_{{\bf N}^d} ^{2}+\mu\left\| \nabla \bm{\phi}^{\bf N}\right\| _{{\bf N}^d}^{2}+(\lambda+\mu)\left\| \mbox{div} \bm{\phi}^{\bf N}\right\| _{{\bf N}^d}^{2}\bigg).
		\end{align}
		This energy is conserved, i.e. $ E^{\bf N}(t)= E^{\bf N}(0)$ for all $t>0$. As usual, this is obtained just multiplying \eqref{adj-dis} by $ \bm{\phi}^{\bf N}_{t} \omega_{\bf i} $, adding in ${\bf i\in I}$ and integrating with respect to time. It is worth noting that the corresponding conservation of energy for the continuous system \eqref{adjcont} requires an integration by parts in space. In order to do the same, here we have to transform the sum in ${\bf i\in I}$ into a ${\bf x}\in \Omega$-integral. This can be done here since the integrand is a polynomial of degree ${\bf 2N-1}$ and the quadrature  formula \eqref{intn} is exact for such polynomials. 
		
		The norm equivalence in \eqref{equiv-n}, together with the conservation of the discrete energy gives,
		\begin{align*}
			\dfrac{1}{2}\int_{0}^{T}\int_{\Omega}\bigg( \left| \bm{\phi}^{\bf N}_{t}\right|^{2}+\mu\left| \nabla \bm{\phi}^{\bf N}\right|^{2}+(\lambda+\mu)\left| \mbox{div} \bm{\phi}^{\bf N}\right| ^{2}\bigg)d{\bf x} dt \geq \frac1{C_2}\int_{0}^{T}E^{\bf N}(t)d{\bf x} dt = \frac{T}{C_2} E^{\bf N}(0).
		\end{align*} 
		We now estimate the terms in the right hand side of \eqref{lions1}. We start with $\left| X+\dfrac{d-1}{2}Y \right|$. This is a quantity evaluated at the times $t=0,T$. This quantity is easily estimated by the discrete energy (using the equivalence of the norm in \eqref{equiv-n}) which is conserved in time. In particular we find,
		\begin{align} \label{lions2}
			\begin{split}
				\left| X+\dfrac{d-1}{2}Y \right|&\le \dfrac{4\sqrt{d}}{\sqrt{\mu}}E^{\bf N}(0).
			\end{split}
		\end{align} 
		Concerning the second term in the right hand side of \eqref{lions1}, we use the fact that on the boundary $ \bm{\nu\cdot}\dfrac{\partial \bm{\phi}^{\bf N}}{\partial {\bm\nu}}= \mbox{div}\bm{\phi}^{\bf N}$, to obtain (see Lions \cite{cit44b}) 
		\begin{align}
			&\int_{0}^{T}\int_{\Gamma} \left( \mu\dfrac{\partial \bm{\phi}^{\bf N}}{\partial {\bm\nu}} +(\lambda+\mu){\bm\nu}\mbox{div} \bm{\phi}^{\bf N}\right)^{2}d\gamma dt
			\ge \mu\int_{0}^{T}\int_{\Gamma}\left( \mu\left|\dfrac{\partial \bm{\phi}^{\bf N}}{\partial {\bm\nu}}\right|^{2}+(\lambda+\mu)\left|\mbox{div}\bm{\phi}^{\bf N}\right|^{2}\right)  d\gamma dt.
		\end{align} 
		Finally, we turn to estimate the last term on the right hand side in \eqref{lions1}. It is enough to consider one of the $2d$ faces of the domain $ \Omega =(-1,1)^d \subset \mathbb{R}^{d}$. We focus on $ \Gamma_1=\{{\bf x}\in \overline{\Omega}, \mbox{ s.t. } x_1=1\}$ . Let us denote $ {\bf  I}_{\Gamma_1}$  the set of indexes corresponding to the collocation nodes on the boundary $ \Gamma_{1}$. For ${\bf i} \in {\bf  I}_{\Gamma_1}$, the Lagrangian basis can be written as $ \Psi_{\bf i}({\bf x})=\Psi_N(x_1)\Psi'_{\bf i}({\bf x'}) $ where $ {\bf x'}=(x_2,...,x_d)\in \mathbb{R}^{d-1} $ and $ \Psi'_{\bf i}({\bf x'})=\prod_{j=2}^{d}\Psi_{k_j}(x_j) $. Therefore, for $\bf x \in \Gamma_1$, ${\bf x}=(1,{\bf x}')$, and we have,
		\begin{align*}
			& \sum_{{\bf i\in I}_{\Gamma_1}} \bigg( \mu\dfrac{\partial^{2} \bm{\phi}^{\bf N} }{\partial \bm{\nu}^{2}}+(\lambda+\mu)\bm{\nu}\dfrac{\partial \mbox{div}\bm{\phi}^{\bf N} }{\partial \bm{\nu}}\bigg)(t,{\bf P_{\bf i}})\Psi_{\bf i}(\bf x) \\
			& \quad = \Psi_N(x_1) \sum_{{\bf i\in I}_{\Gamma_1}} \left( \mu\dfrac{\partial^{2} \bm{\phi}^{\bf N} }{\partial {\bm \nu}^{2}}+(\lambda+\mu){\bm \nu}\dfrac{\partial \mbox{div}\bm{\phi}^{\bf N} }{\partial {\bm \nu}}\right) (t,{\bf P_{\bf i}})\Psi_{\bf i}'(\bf x')
			\\
			& \quad =  \Psi_N(x_1)  \left( \mu\dfrac{\partial^{2} \bm{\phi}^{\bf N} }{\partial {\bm \nu}^{2}}+(\lambda+\mu){\bm \nu}\dfrac{\partial \mbox{div}\bm{\phi}^{\bf N} }{\partial {\bm \nu}}\right)(t,1,{\bf x'}) .
		\end{align*}		
		We now replace this in the last term on the right hand side in \eqref{lions1}. Using the Young's inequality,
		\begin{align}\label{estGamma12}
			\begin{split}
				&\left| \int_{0}^{T}\int_{\Omega}\left( \left( \mu\dfrac{\partial^{2} \bm{\phi}^{\bf N} }{\partial {\bm \nu}^{2}}+(\lambda+\mu){\bm \nu}\dfrac{\partial \mbox{div}\bm{\phi}^{\bf N} }{\partial {\bm \nu}}\right)(t,1,{\bf x'})\Psi_{N}(x_1)\right)\bm{\cdot} \left(  m_j\dfrac{\partial \bm{\phi}^{\bf N}}{\partial x_{j} }\right)({\bf x}) d{\bf x} dt \right| \\
				&\le \int_{0}^{T}C_\varepsilon \int_{-1}^{1}\int_{[-1,1]^{d-1}}\left| \left( \mu\dfrac{\partial^{2} \bm{\phi}^{\bf N} }{\partial {\bm \nu}^{2}}+(\lambda+\mu){\bm \nu}\dfrac{\partial \mbox{div}\bm{\phi}^{\bf N} }{\partial {\bm \nu}}\right) (t,1,{\bf x}')\right|^2\left| \Psi_{N}(x_1)\right| ^2 dx_1d{\bf x}' dt \\
				& + \int_{0}^{T}\varepsilon \int_{\Omega}\bigg| m_j\dfrac{\partial \bm{\phi}^{\bf N}}{\partial x_{j} }\bigg|^2d{\bf x} dt \le C_\varepsilon |\Psi_{N}|^2_{L^2(-1,1)} \int_{0}^{T}\int_{\Gamma_1}\left| \mu\dfrac{\partial^{2} \bm{\phi}^{\bf N} }{\partial {\bm \nu}^{2}}+(\lambda+\mu){\bm \nu}\dfrac{\partial \mbox{div}\bm{\phi}^{\bf N} }{\partial {\bm \nu}}\right|^2d\gamma dt \\
				&+ \varepsilon (\mbox{sup}_{{\bf x}\in \Gamma}m_j\nu_j)^2\int_{0}^{T}\int_{\Omega}| \nabla \bm{\phi}^{\bf N}|^2d{\bf x} dt.
			\end{split}
		\end{align}
		Taking into account the norm equivalence in \eqref{equiv-n}, $\left|\Psi_{N} \right|^2 _{L^2(-1,1)}\le \left\| \Psi_{N} \right\|^2 _N= \omega_{N} $, the conservation of the discrete energy proved above  in \eqref{disenergy}, the fact that $sup_{{\bf x}\in \Gamma}\ { m_j }{\nu_j}=2\sqrt{d}$ and \eqref{estGamma12}, we obtain
		\begin{align}\label{est-boundary}
			\begin{split}
				&\left|\int_{0}^{T}\int_{\Omega}\left(  \sum_{{\bf i\in I}_{\Gamma_1}}\left( \mu\dfrac{\partial^{2} \bm{\phi}^{\bf N} }{\partial {\bm \nu}^{2}}+(\lambda+\mu){\bm \nu}\dfrac{\partial \mbox{div}\bm{\phi}^{\bf N} }{\partial {\bm \nu}}\right)(t,{\bf P_i})\Psi_{\bf i}({\bf x})\right) \bm \cdot \left(  m_j\dfrac{\partial \bm{\phi}^{\bf N}}{\partial x_{j} }\right) ({\bf x}) d{\bf x} dt\right|\\
				&\le C_\varepsilon\omega_{N}\int_{0}^{T}\int_{\Gamma_{1}}\left|\mu\dfrac{\partial^{2} \bm{\phi}^{\bf N} }{\partial {\bm \nu}^{2}}+(\lambda+\mu){\bm \nu}\dfrac{\partial \mbox{div}\bm{\phi}^{\bf N} }{\partial {\bm \nu}} \right|^2d\gamma dt+\dfrac{8d}{\mu}\varepsilon TE^{\bf N}(0).
			\end{split}
		\end{align}		
		An analogous estimate holds for the other $ 2d-1 $ terms of the boundary in the right hand side of \eqref{lions1}.
		It follows from \eqref{lions1}-\eqref{est-boundary} and the fact that $ \omega_{N}=\omega_{0}=\dfrac{2}{N(N+1)} $ (see  \cite{ex5}, Chapter 2), 
		\begin{align}\label{final-est}
			\begin{split}
				&\bigg(\dfrac{T}{C_{2}}-\dfrac{4\sqrt{d}}{\sqrt{\mu}}-\dfrac{16d^{2}\varepsilon T}{\mu}\bigg)E^{\bf N}(0)\le  \dfrac{\sqrt{d}}{\mu}\int_{0}^{T}\int_{\Gamma} \left| \mu \dfrac{\partial \bm{\phi}^{\bf N}}{\partial {\bm \nu}} +(\lambda+\mu){\bm \nu}\mbox{div}\bm{\phi}^{\bf N} \right| ^{2}d\gamma dt\\
				&+\dfrac{4dC_\varepsilon}{N^2}  \int_{0}^{T}\int_{\partial\Omega}\left| \mu \dfrac{\partial^{2} \bm{\phi}^{\bf N}}{\partial {\bm \nu}^{2}} +(\lambda+\mu){\bm \nu}\dfrac{\partial \mbox{div} \bm{\phi}^{\bf N}}{\partial {\bm \nu}}\right|^{2}d\gamma dt.
			\end{split}
		\end{align} 
		The fact that we can replace the discrete energy $E^{\bf N}(0)$, by the $(H^1_0)^d \times (L^2)^d$ norm of the initial data is a consequence of the equivalence of the discrete $L^2$-norm  in \eqref{equiv-n} for polynomials of degree ${\bf N}$.  		
		Inequality \eqref{obser-inq} holds as long as  $\dfrac{T}{C_2}-\dfrac{16d^{2}T\varepsilon}{\mu}-\dfrac{4\sqrt{d}}{\sqrt{\mu}}>0$. As $\varepsilon $ can be chosen arbitrarily small and $C_2=(2+N^{-1})^d$ we have the condition $T>\dfrac{4\sqrt{d}(2+N^{-1})^d}{\sqrt{\mu}}$.
	\end{proof}
	
		\section{Application: boundary control}
	
	\hspace*{0.5cm}Let us consider the following boundary control problem for the linear elasticity: Given the initial data $ ({\bf u^0,u^1})\in (L^2(\Omega))^d\times(H^{-1}(\Omega))^d $ and $T> 4\sqrt{d}/\sqrt{\mu}$, find a control $ {\bf f}\in(L^2(0,T;\Gamma) )^d $ such that the solution $ {\bf u}\in (C([0,T]; H^{1}_0(\Omega)))^d\cap  (C^1([0,T]; L^{2}(\Omega)))^d  $ of the system  
	\begin{align}\label{eq1}
		\begin{cases}
			\left({\bf u }_{tt}-\Delta^* {\bf u }\right)(t,{\bf x})={\bf 0},& (t,{\bf x}) \in (0,T)\times\Omega, \\
			{\bf u }(t,{\bf x})={\bf f }(t,{\bf x})\chi_{\Gamma }({\bf x}), & (t,{\bf x})\in (0,T)  \times \partial\Omega, \\
			({\bf u}(0,{\bf x}), {\bf u}_{t}(0,{\bf x}))=(\bf {u^0(x),u^1(x)}),\hspace*{0.5cm}& {\bf x} \in  \Omega,
		\end{cases}
	\end{align}
	satisfies the null controllability condition
	$ {\bf u}(T,{\bf x})={\bf u}_{t}(T,{\bf x})={\bf 0}, \ {\bf x }\in\Omega$. Here $ \chi_{\Gamma } $ is the characteristic function of the set $ \Gamma \subset \partial \Omega $. It is well-known that a control ${\bf f}$ exists (\cite{cit6} and \cite{cit44b}). Moreover, among all the controls, the one with minimal $L^2-$norm is unique. Our main objective is to approximate this control ${\bf f}$. 
	Let us denote $ \partial\Omega=\cup_j^{2d} \Gamma_j $ where $ \Gamma_j=\{{\bf x}\in \overline{\Omega}, \mbox{ s.t. } x_j=1\}$ and $ \Gamma_{d+j}=\{{\bf x}\in \overline{\Omega}, \mbox{ s.t. } x_j=-1\},\;j=1,...d$. The set $ \mathbb{P}_{\bf N}(\Gamma)$ (resp $ \mathbb{P}_{\bf N}(\Gamma_j),j=1,...,2d$)  is the restriction of $  \mathbb{P}_{\bf N}(\Omega) $  to $ \Gamma $ (resp to $\Gamma_j,j=1,...,2d$). We introduce the following discrete control problem: Given $ {\bf{u^{0,N}}},\; {\bf{u^{1,N}}}$ $ \in (\mathbb{P}^{Di}_{\bf N}(\Omega) )^d $ and $T>\dfrac{4\sqrt{d}(2+N^{-1})^d}{\sqrt{\mu}}$, find $ {\bf f}^{\bf N}\in (L^2(0,T;\mathbb{P}_{\bf N}(\Gamma)) )^d$,\;$ {\bf g}^{\bf N}_j \in (L^2(0,T;\mathbb{P}_{\bf N}(\Gamma_j)))^d,j=1,...,2d$, such that the solution  $ {\bf u^{N}}\in (C^1 ([0,T]; \mathbb{P}_{\bf N}(\Omega)))^d  $ of system
	\begin{align}\label{control-d}
		\begin{split}
			\begin{cases}
				\left( {\bf u}_{tt}^{\bf N}-\Delta^* {\bf u^{N}}\right) (t,{\bf  P_{\bf i}})= \sum_{j=1}^{2d} G_{j}^{\bf N}(t,{\bf P_i}), 
				\ & (t,{\bf  P_{\bf i}})\in(0,T) \times C^{\Omega},\\
				{\bf u^{N}}(t,{\bf  P_{\bf i}})={\bf f^N}(t,{\bf P_i})\; \chi_{\Gamma }({\bf P_i}),\ &  (t,{\bf P_i })\in(0,T) \times C^{\partial\Omega}, \\
				{\bf u^{N}}(0,{\bf P_i })={\bf u^{0,N}}({\bf P_i }),\quad {\bf u}_{t}^{\bf  N}(0,{\bf P_i })={\bf u^{1,N}}({\bf P_i }),	\ & {\bf  P_{\bf i}} \in C^{\Omega},
			\end{cases}
		\end{split}
	\end{align}
	satisfies $ {\bf u^{N}}(T,{\bf  P_{\bf i}})={\bf u}_{t}^{\bf N}(T,{\bf P_ i})={\bf 0},$ at ${\bf P_i} \in C^{\Omega}.$  
	Here $ G_{j}^{\bf N}\in (C^1([0,T];\mathbb{P}_{\bf N}(\Omega)))^d$ depend on $ g_j^{\bf N},\;j=1,...,2d $. For example, consider $ j=1 $ and $ j=d+1 $, which corresponds to $ \Gamma_1,\Gamma_{d+1} $ respectively. 
	For any point	 
	$ {\bf P_i}=(x_1^{k_1},x_2^{k_2},...,x_d^{k_d})\in C^{\Omega} $, we write the proyections of this point on these two opposite sides of the domain, 
	$ {\bf P}^1_{\bf i}=(1,x_2^{k_2},...,x_d^{k_d}) \in \Gamma_1$, $  {\bf P}^{2}_{\bf i}=(-1,x_2^{k_2},...,x_d^{k_d}) \in \Gamma_{d+1}$,  and 
	\begin{align}\label{condc}
		\begin{split}
			\begin{cases}
				G_{1}^{\bf N}(t,{\bf  P_{\bf i}})= A_1 \tilde{h}_1(x_1^{k_1}) g_1^{\bf N}(t,{\bf P}^1_{\bf i}),\qquad	G_{d+1}^{\bf N}(t,{\bf  P_{\bf i}})= A_1 \tilde{h}_2 (x_1^{k_1}) g_{d+1}^{\bf N}(t,{\bf P}^{2}_{\bf i})\\
				A_j\in \mathcal{M}_{d\times d}\; \mbox{diagonal with components }\;a_{kk}=\mu+(\mu+\lambda)\delta_{kj}\\
				\tilde{h}_1(s)=\dfrac{1}{\sqrt{ \omega_N}}\left(h^{1}_{ss}+\dfrac{{\Psi_{N,s}}}{ \omega_{N}}\right)(s),\qquad 	\tilde{h}_2(s)=\dfrac{1}{\sqrt{ \omega_0}}\left(h^{2}_{ss}-\dfrac{{\Psi_{0,s}}}{ \omega_{0}}\right)(s)\\
				h^{1},h^{2}\in \mathbb{P}_{ N}^{Di}(-1,1),\quad h^{1}(s^k)=\dfrac{1+s^k}{2}, \quad  h^{2}(s^k)=\dfrac{1-s^k}{2},\quad s^k \in C^{(-1,1)},
			\end{cases}
		\end{split}
	\end{align}
	where $ C^{(-1,1)}  $ are the interior LGL nodes in $ (-1,1) $. Analogous formulas define the other $G_j^{\bf N}$ in terms of the controls ${\bf g}_j^{\bf N}$ for $j=2,...,2d$. 
	
	
	Comparing the control problems \eqref{eq1} and \eqref{control-d} we see that the discrete version has $ 2d $ extra controls ${\bf g}^{\bf N}_j$ depending on $\Gamma_j$, $j=1,...,2d$, that together can be interpreted as a single control in the whole boundary $ {\bf g}^{\bf N} \in (L^2(0,T;\mathbb{P}_{\bf N}(\partial\Omega)))^d$. This 'numerical' boundary control is associated with the last term in the inequality \eqref{obser-inq} and it is necessary to obtain a bounded sequence of discrete controls ${\bf f^N}$ as ${\bf N}\to \infty$. In fact, the existence of discrete controls and their uniform bound (with respect to $N$) in terms of the initial data is a direct consequence of the observability result established in Theorem \ref{obser} above.
	We refer to \cite{Som} where this is detailed in the case of the wave equation. For this much simpler wave equation the authors go further and prove that, under some technical hypotheses, $ {\bf f^N\to f} $ in $(L^2(0,T;\Gamma ))^d$ and the numerical control $ {\bf g^N} $ vanishes as $ {\bf N \to \bm \infty }$. For the elasticity system such convergence  result is more difficult to prove. In fact, if we try to apply the general theory established in \cite{EZ} to recover the convergence of controls, some important hypotheses related with the convergence of the solutions of \eqref{adjcont} must be checked. In particular, it should be established that for a convergent sequence of discrete initial data
	$ (\bm {\phi^{0,N},\phi^{1,N}})\to  (\bm {\phi^{0},\phi^{1}})$ in $ (H_0^1(\Omega))^d\times (L^2(\Omega))^d $ the right hand side in \eqref{obser-inq} converges to the right hand side in \eqref{obs_ineq}. More precisely, 
	\begin{eqnarray*}
		&&\mu \dfrac{\partial\bm{\phi}^{\bf N}}{\partial{\bm \nu}}+(\lambda+\mu){\bm \nu} \mbox{div} \bm{\phi}^{\bf N}\to \mu \dfrac{\partial\bm{\phi}}{\partial{\bm \nu}}+(\lambda+\mu){\bm \nu} \mbox{div} \bm{\phi}, \quad \mbox{ in } L^2(0,T;\Gamma),\\
		&&\mu \dfrac{\partial^{2}\bm{\phi}^{\bf N}}{\partial{\bm \nu}^{2}}+(\lambda+\mu){\bm \nu}\dfrac{\partial\mbox{div}{\bm\phi}^{\bf N}}{\partial{\bm \nu}} \to {\bf 0}, \quad \mbox{ in } L^2(0,T;\Gamma).
	\end{eqnarray*}
	Such convergent results are not standard with the usual numerical analysis techniques and require further investigation.      
	
	Here we present numerical evidences of the convergence of the discrete control to the continuous one.	We consider the following initial conditions,  
	${\bf u^0}(x_{1}, x_{2}) = (0.2 \sin (\pi (x_{1}+1)/2) \sin (\pi (x_{2}+1)/2), 0.2 \sin (\pi (x_{1}+1)/2) \sin ( \pi (x_{2}+1)/2)), $ ${\bf u^1}(x_{1}, x_{2}) = (0, 0).$
	We take $ \lambda= 0.5 $, $\mu = 4$ and final time $T = 3$ with time step $ \Delta t= 0.01$. Note that the time control is only slightly greater than the minimal control time for the continuous elasticity system $ T> \dfrac{4\sqrt{2}}{\sqrt{\mu}} $ (see \cite{cit44b}) and lower than the time given by the uniform discrete observability inequality in Theorem \ref{obser} (which is probably not optimal).
	
	\begin{minipage}{\textwidth}
		\begin{minipage}{0.45\textwidth}
			\centering
			\includegraphics[width=70mm]{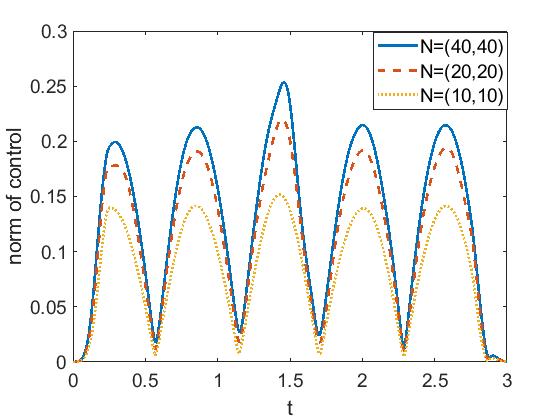}
			\captionof{figure}{Time behavior of $| f^N|_{L^2(\Gamma)}$.}
			\label{fig_1}
		\end{minipage}
		\hfill
		\begin{minipage}{0.55\textwidth}
			\centering
			\begin{tabular}{ |c|c|c| } 
				\hline
				$ {\bf N} $ & $\left| {\bf f^N}\right| _{(L^{2}(\Gamma))^2} $& $\left| {\bf g}^{\bf N}\right| _{(L^{2}(\partial\Omega))^2} $\\
				\hline
				$ (10,10) $	&$ 1.6\times10^{-1}  $& $ 4.8\times 10^{-3}  $\\ 
				$ (20,20) $ & $ 2.2\times10^{-1}  $ &  $ 2.5\times 10^{-3}  $\\
				$ (40,40) $&  $2.5\times10^{-1}$   & $ 8.2\times 10^{-4}$ \\ 
				\hline
			\end{tabular}
			\captionof{table}{Norm of the discrete controls.}
			\label{table_1}
		\end{minipage}
	\end{minipage}
	
	In Figure \ref{fig_1} we show the behavior of the norm of control $ {\bf f^N} $ in time for different values of $ {\bf N} $. Table \ref{table_1} illustrates the behavior of the norm for the controls when the degree of the polynomials $ \bf N $ grows. We observe that the boundary control ${\bf f^N}$ remains bounded and should converge to a continuous control ${\bf f}$ while the numerical artificial control ${\bf g}^{\bf N}$  vanishes as ${\bf N} \to \infty$.

\section*{Declarations}

{\bf Funding:} The authors were supported by grant PID2021-124195NB-C31 from the Spanish Governement (MICINN).

\bigskip

\noindent {\bf Conflict of interest:} The authors declare no conflict of interest. 

\bibliography{Reference2.bib}
\end{document}